\documentclass[11pt,letterpaper,leqno,oneside,fleqn]{article}

\input SDD.sty
\def\version{18 April 2014}

\begin{document}

\Heading
Remarks on Kneip's linear smoothers
\endHeading

\begin{center}
S\"oren R. K\"unzel\\
David Pollard\\
Dana Yang\\
\ \\
Statistics Department, Yale University\\
\version
\end{center}
\TOC

\bigskip

\section{Introduction}\label{intro}
We have been trying to understand the analysis provided by \citet{Kneip1994AnnStat}.
In particular we want to persuade ourselves that his results imply the oracle inequality stated 
by \citet[Lecture~8]{Tsybakov2014yale}. 

This note contains our reworking of Kneip's ideas.
We refer to page~x of Kneip's paper as Kx.
For $n\times n$ symmetric matrices we write $A \preccurlyeq B$ to mean that $B-A$ is \psd/. 
Also we write~$|\cdot|$ for the usual Euclidean length in~$\RR^n$, that is,~$|x|^2=\SUM_{i\le n}x_i^2$.

Following Kneip, we consider an observed  $n\times 1$ random vector $y=\mu+\xi$ with unknown $\mu$ 
and error~$\xi$
(with independent components)
with $\PP \xi=0$ and $\var(\xi)=\sig^2 I_n$.  
We assume that $\xi\sim N(0,\sig^2I_n)$.
Kneip(K844, statement of Theorem~1) assumed subgaussianity.
The possible estimators are of the form~$Sy$, with~$S$
in a specified set $\ss$ of $n\times n$ (symmetric) \psd/ smoothing matrices
 that is totally ordered under the semi-definite ordering~$\preccurlyeq$, with $0\preccurlyeq S \preccurlyeq I_n$ for all $S\in\ss$.

Kneip considered the estimator  $\Shat y$  with
\beqN
\Shat = \ARGMIN_{S\in\ss}\Ghat(S)
\qt{where } \Ghat(S) = |y-S y|^2 + 2\sig^2 \trace(S).
\eeqN
Here and subsequently we omit multiplicative factors of~$n^{-1}$ that Kneip used. This selection procedure is the well known \textit{Mallows' $C_p$}.

The analysis and the statement of Kneip's main result involve two related processes, which we define for all positive semi-definite matrices~$S$:
\bAlign
G_\mu(S) &:= |\mu - S y|^2 
\\
M_\mu(S) &:= \PP G_\mu(S) = |\mu-S\mu|^2 + \sig^2\trace(S^2).
\eAlign
Following Kneip, we assume that the minimium of $M_\mu$ over the set~$\ss$ is achieved at the matrix~$\Smu$ in~$\ss$ and define 
\beqN
m^*=M_\mu(S_\mu) = \MIN_{S\in\ss}M_\mu(S)
\eeqN
We ignore all questions of whether mininima are achieved and whether~$\Shat$ is measurable.

\begin{theorem} 	\label{K844} (K844)
There exist constants $C_1$ and $C_2$ that depend only on $\sig^2$ for which for all $\mu$ in $\RR^n$,
\beqN
\PP\{ |G_\mu(\Shat) -  G_\mu(\Smu)| \ge \max\left( x^2,x\rootmstar\right)\}\le C_1 e^{-C_2 x}
\qt{for }x\ge 0.
\eeqN
\end{theorem}


\begin{corollary} 	\label{oracle}
There exist constants $C_3$ and $C_4$ that depend only on $\sigma^2$ for which for all $\mu$ in $\RR^n$,
\beqN
\PP G_\mu(\Shat)  \le m^* +C_3 \rootmstar + C_4.
\eeqN
\end{corollary}

\startDP
The Corollary is equivalent to 
\endDP
\beqN
\PP G_\mu(\Shat)  \le (1+\eps)m^*+ C_0/\eps +C_4
\qt{for all $\eps>0$ and $C_0=C_3^2/4$},
\eeqN
\startDP a minor modification of the  oracle inequality stated by \citet[Lecture~8]{Tsybakov2014yale}. For~$\eps$ in a bounded range the~$C_4$ can be absorbed into the previous term.\endDP

The proof of the Theorem makes extensive use of the properties of the metric~$d$ defined 
\startDP
on the set of
\endDP all positive semi-definite matrices~$S_1$ and~$S_2$ by
\beqN
d^2(S_1,S_2) = \PP|S_1y-S_2y|^2 = |(S_1-S_2)\mu|^2 + \sig^2\trace(S_1-S_2)^2
\eeqN
(Note that $d^2(S_1,S_2)=nq_\mu^2(S_1,S_2)$
for the~$q_\mu$ defined near the bottom of~K842.) 
In particular, the proof relies crucially on a bound 
\startDP
(see Section~\ref{packing})
for the packing numbers of subsets of~$\barSS$, a set of \psd/ 
matrices
\endDP
 that contains~$\ss$ as a subset. The arguments rely on the total ordering of~$\ss$ to parametrize~$\ss$  by a
subset of the real line.

\section{Outline of the Proofs}\label{top.level}
\startDP
To prove Theorem~\cref{K844} we first show that
\bAlign
\Ghat(S) &\approx M_\mu(S) + \text{term not depending on $S$}
\\
G_\mu(S) &\approx M_\mu(S) + \text{term not depending on $S$}.
\eAlign
More precisely, with
\bAlign
D_\mu(S) &:= G_\mu(S) - M_\mu(S) 
\AND/
\Dhat(S) := \Ghat(S)-M_\mu(S),
\eAlign
we show: 
There exist positive constants $C_1$, $C_2$, depending only on~$\sig^2$ for which, for every $r>0$,
\bAlignL
\PP\{\exists S\in\ss:& 
|\Dhat(S)-\Dhat(\Smu)|> L(S,x,r)\} 
\le C_1e^{-C_2 x},
\label{Dhat.max}
\\
\PP\{\exists S\in\ss:& 
|D_\mu(S)-D_\mu(\Smu)|> L(S,x,r)\} 
\le C_1e^{-C_2 x},
\label{Dmu.max}
\\
\label{L.def}
&\qt{where }L(S,x,r) = \left[d^2(S,\Smu)+r^2\right]x/r .
\eAlignL
The proof of these inequalities (in  Section~\ref{technical}) uses a chaining argument based on control of the increments of both the~$\Dhat$ and~$\Dmu$ processes, together with a bound on the packing numbers that derives from the total ordering of~$\ss$.

We also make use of an inequality (cf. K843, Proposition~1)
related to
the growth of~$M_\mu(S)-M_\mu(\Smu)$ as~$d(S,\Smu)$ increases. For that we need the matrix analog of the inequality $\a^2+\beta^2\ge(\a-\beta)^2$ for nonnegative real numbers.
\endDP

\begin{lemma}  \label{square.diff}
If $S_1$ and $S_2$ are symmetric, positive semi-definite matrices that commute then
$(S_1-S_2)^2 \preccurlyeq S_1^2+S_2^2$.
\end{lemma}
\Proof
We want to show that the matrix
$$(S_1^2+S_2^2) - (S_1-S_2)^2=2S_1S_2$$
is positive semi-definite. Let $U$ be an orthogonal matrix that simultaneously diagonalizes $S_1$ and $S_2$ to $\Lambda_1$ and $\Lambda_2$. Then for any vector $\alpha$ in $\RR^n$, we have
\beqN
\alpha'S_1S_2\alpha=(U\alpha)'\Lambda_1\Lambda_2(U\alpha),
\eeqN
which is nonnegative because 
the elements of the diagonal matrix~$\Lam_1\Lam_2$ are all nonnegative.
\endProof
As a direct consequence of the Lemma,
\bAlignL
M_\mu(S_1)&+M_\mu(S_2) 
\notag\\
&= \mu'\left[(I_n-S_1)^2+(I_n-S_2)^2\right]\mu + \sig^2\trace\left[S_1^2+S_2^2\right]
\notag\\
&\ge \mu'(S_1-S_2)^2\mu + \sig^2\trace(S_1-S_2)^2 = d^2(S_1,S_2).
\label{M.incr}
\eAlignL
In particular, if $d^2 (S,\Smu) \ge 3m^*$ then $d^2 (S,\Smu) \ge d^2 (S,\Smu)/3+2m^*$, so that~\cref{M.incr} implies
\beq\label{M.growth}
M_\mu(S) - m^* \ge \tfrac13 d^2(S,\Smu)\{d(S,\Smu) \ge \Sqrt{3m^*}\}.
\eeq

\startDY
\Proof 
\startDP(of Theorem \cref{K844})
With $L$ as defined in~\cref{L.def}, define 
\beqN
L(S,x):= L(S,x,r_x)
\qt{where $r_x=\max(\Sqrt{3m^*},7x)$.}
\eeqN
By inequalities~\cref{Dhat.max}
\endDP
 and~\cref{Dmu.max}, we can find a set~$\Om_x$ with probability at least $1-2C_1e^{-C_2 x}$, on which we have
\beq\label{Omx}
\max\left(|\Dhat(S)-\Dhat(\Smu)|,|\Dmu(S)-\Dmu(\Smu)|\right)
\le L(S,x)
\qt{for all $S\in\ss$}.
\eeq
\startDP 
The rest of the proof is just a  deterministic argument  on the set~$\Om_x$.

Define $\dhat=d(\hat{S},S_\mu)$.
Then
\bAlign
\tfrac17(\dhat^2+r_x^2) &\ge 
L(\Shat,x) 
\qt{because $x/r_x<1/7$}
\\
&\ge \Dhat(\Smu) -\Dhat(\Shat)
\qt{by~\cref{Omx}}
\\
&= 
\Ghat(\Smu) - \Ghat(\Shat) + M_\mu(\Shat)-M_\mu(\Smu)
\\
&\ge M_\mu(\Shat)-m^*
\qt{because $\Shat$ minimizes $\Ghat$}
\\
&\ge 
\tfrac13 \dhat^2\{\dhat \ge \sqrt{3m^*}\}
\qt{by~\cref{M.growth}.}
\eAlign
If $\dhat$ were larger than~$r_x$ the last inequality would give
$
\tfrac27 \dhat^2 \ge \tfrac13\dhat^2
$, 
which clearly cannot be true. Thus $\dhat<r_x$ on~$\Om_x$, implying
\beqN
2r_x^2 x/r_x \ge L(\Shat,x) \ge M_\mu(\Shat)-m^* .
\eeqN
In summary,
\beq\label{compared}
\dhat := d(\Shat,\Smu) < r_x
\AND/
M_\mu(\Shat) \le m^* + 2xr_x
\qt{on $\Om_x$.}
\eeq

Combine this inequality with the bound for~$|\Dmu(S)-\Dmu(\Smu)|$ from~\cref{Omx} to deduce that, again on~$\Om_x$,
\bAlign
|G_\mu(\hat{S})-G_\mu(S_\mu)|
&\leq\left(M_\mu(\hat{S})-M_\mu(S_\mu)\right)+|D_\mu(\hat{S})-D_\mu(S_\mu)|
\\
&\le 2xr_x + L(\Shat,x)
\\
&\le  4xr_x.
\eAlign
Thus
\beqN
\PP\{|G_\mu(\hat{S})-G_\mu(S_\mu)|>4xr_x\}\le \PP \Om_x \le 2k_1 e^{-k_2 x}.
\eeqN
This inequality is not quite the result announced in  Theorem~\cref{K844}. However,
\beqN
4xr_x = 4x \max(\sqrt{3m^*},7x)\ge \max(4\sqrt{3},28)\max(x\rootmstar,x^2)
\eeqN
so we get the announced result, for $Z=|G_\mu(\Shat) -  G_\mu(\Smu)|$:
\beq\label{Z.tail}
\PP\{ Z \ge \max( x^2,x\rootmstar\,)\}\le C_1 e^{-C_2 x}
\qt{for }x\ge 0.
\eeq
by adjusting the constants. 
\endDP
%
%
%
\endProof
  
The oracle inequality stated 
\startDP as Corollary~\cref{oracle} is an integrated version of
the tail bound from Theorem~\cref{K844}.\endDP
\endDY
\startSK
\Proof
\startDP
From inequality~\cref{Z.tail} we have $\PP\{Z \ge f(x)\}\le C_1e^{-C_2 x}$ for $x\ge0$, where $f(x) = \max(x^2,x\rootmstar)$, which gives
\bAlign
|\PP G_\mu(\Shat) -m^*| &\le \PP Z =  \int_0^\infty \PP\{Z>t\}\,dt
=\int_0^\infty \PP\{Z\ge f(x)\}f'(x)\,dx
\\
&\le \int_0^{\infty} 
 \max(2x, \sqrt{m^*})C_1e^{-C_2 x}\,dx
\\
&\le \frac{C_1}{C_2}\sqrt{m^*} +2 \frac{C_1}{C_2^2}e^{-C_2\sqrt{m^*}} + \frac{C_1}{C_2}\sqrt{m^*} e^{-C_2\sqrt{m^*}}
\\
&\le C_3\rootmstar + C_4
\eAlign
\endDP
for new constants $C_3=2C_1/C_2$ and $C_4=2C_1/C_2^2$.

\endProof
\endSK

\section{Technical Stuff}\label{technical}
\startDP
This section proves the inequalities~\cref{Dhat.max} and~\cref{Dmu.max},
\bAlign
\PP\{\exists S\in\ss:& 
|\Dhat(S)-\Dhat(\Smu)|> L(S,x,r)\} 
\le C_1e^{-C_2 x},
\\
\PP\{\exists S\in\ss:& 
|D_\mu(S)-D_\mu(\Smu)|> L(S,x,r)\} 
\le C_1e^{-C_2 x},
\\
&\qt{where }L(S,x,r) = \left[d^2(S,\Smu)+r^2\right]x/r ,
\eAlign
by means of a chaining argument with stratification. The necessary ingredients are the control of increments of the $\Dhat$ and $\Dmu$ processes and bounds on packing numbers.

\endDP
	\subsection{Exponential bounds for increments}\label{increments}
\startDP
The next Lemma is all we need to control the increments
of the $\Dhat$ and~$\Dmu$ processes under the assumption of gaussian errors.
First we expand each process into sums of simpler processes.
\bAlign
D_\mu(S)= & G_\mu(S)-M_\mu(S) = X_1(S)+X_2(S)\\
\hat{D}(S)= & \hat{G}(S)-M_\mu(S) = X_3(S)+X_4(S)+n\sigma^2. 
\eAlign
where
\beq\label{Xi.defs}
\begin{split}
X_1(S)  &= \xi'S^2\xi-\sigma^2\trace(S^2)
\\
X_2(S) &= -2\mu'(S-S^2)\xi
\\
X_3(S) &= \xi'(I-S)^2\xi-\sigma^2\trace(I-S)^2
\\
X_4(S) &= 2\mu'(I-S)^2\xi
\end{split}
\eeq
Notice that each $X_i(S)$ is either a linear or quadratic function of $\xi$. 
\endDP

\begin{lemma}  \label{exp.bnd} (Compare with Kneip's Lemma~2, K852)
Suppose $z\sim N(0,I_n)$. For each vector of constants~$a$ and each symmetric matrix~$A$,
\bAlign
\PP\{z'a \ge w|a|\} &\le \exp(-w^2/2)\\
\PP\{z'Az - \trace(A) \ge w \Sqrt{\trace(A^2)}\,\} &\le 2e^{-w/4}
\eAlign
for each $w\ge 0$.
\end{lemma}
\Proof
The first inequality is just the usual bound for $N(0,1)$ tails.
(It extends easily to the subgaussian case.)
For the second inequality write $A$ as $L'\diag(\lam_1,\dots,\lam_n)L$, 
with~$L$ orthogonal. Write~$\kappa$ for~$\Sqrt{\trace(A^2)}=|\lam|$. Then~$x=Lz\sim N(0,I_n)$. With $t=1/(4\kappa)$,
\bAlign
\PP\{z'Az - \trace(A) &\ge w \Sqrt{\trace(A^2)}\}
\\
&= \PP \{\SUM_i \lam_i(x_i^2-1) \ge w\kappa\}
\\
&\le e^{-tw\kappa}\PROD_i \PP\exp\left( -t\lam_i +t\lam_i x_i^2\right)
\\
&= e^{-w/4}\exp\left( \SUM_i\left(-t\lam_i -\tfrac12\log(1-2t\lam_i)\right)\right).
\eAlign
As $\max_i|2t\lam_i|\le 1/2$, we have
$
-\tfrac12\log(1-2t\lam_i) \le t\lam_i +\tfrac12(2t\lam_i)^2,
$
which leaves $2t^2\SUM_i\lam_i^2= 1/8<\log 2$ in the exponent. 
\endProof

\startDP
The argument for the quadratic form comes from \citet[Lemma~3]{NolanPollard87Uproc1}. For subgaussian errors Kneip calculated moments,  resulting in a bound  similar to an earlier result of
\citet{HansonWright1971AMS}. The  \citet{RudelsonVershynin2013Arxiv} method provides a simpler derivation.

\endDP

The $\exp(-w^2/2)$ bound 
\startDP
for~$z'a$
\endDP
 is more than we need. The inequality
\beqN
\min\left( 1,2e^{-w^2/2}\right)\le 4\exp(-w/4)\qt{for all }w\ge0.
\eeqN
shows that all the increments of the $X_i$ processes from~\cref{Xi.defs} satisfy inequalities of the form
\beq\label{Xi.incr}
\PP\{|X_i(S_1)-X_i(S_2)| >  d(\theta_1,\theta_2)x\}
\le C_1e^{-C_2 x}
\qt{for }x\ge 0,
\eeq
for  constants $C_1$ and~$C_2$.


	\subsection{Packing bounds}\label{packing}
The assumption on~$\ss$ ensures the matrices can be diagonalized by a fixed rotation: $S=U'\Lam(S)U$ with $U$ orthogonal and
\beqN
\Lam(S)= \diag\left( \lam_1(S),\dots, \lam_n(S)\right)
\eeqN
The total ordering ensures that each $S\in\ss$ is uniquely determined by its trace. The set~$\ss$ can be parametrized as $S_\theta$, with $\theta\in\Theta
\subseteq [0,n]$,  where
\beqN
\Lam(S_\theta) = \Lam(\theta) =\diag\left( \lam_1(\theta),\dots, \lam_n(\theta)\right)
\qt{and $\theta = \SUM_{i\le n}\lam_i(\theta)$}.
\eeqN
The maps $\theta\mapsto \lam_i(\theta)$ are increasing, for each~$i$.
As Kneip showed (by interpolation, K857), $\ss$ can be embedded into a larger family of \psd/ matrices
$\barSS=\{S_\theta:\theta\in\Thetabar\}$ with $S_\theta=U'\Lam(\theta)U$ 
and~$\theta\mapsto \lam_i(\theta)$ continuous and nondecreasing from~$\Thetabar=[0,n]$ onto~$[0,1]$. The monotonicity of $\theta\mapsto \lam_i(\theta)$  simplifies calculation of packing/covering numbers for subsets of~$\Thetabar=[0,n]$ under the metric~$d$.
 Recall that 
\beqN
d^2(\theta_1,\theta_2) =\SUM_i (\rho_i^2+\sig^2)|\lam_i(\theta_1)-\lam_i(\theta_2)|^2
\eeqN
and $\theta\mapsto \lam_i(\theta)$ is nondecreasing. If $a\le t_1<t_2<\dots<t_N\le b$ then
\beqN
|\lam_i(b)-\lam_i(a)|^2 \ge \left(\SUM_{j=2}^N \lam_i(t_j)-\lam_i(t_{j-1}) \right)^2 \ge \SUM_j\left| \lam_i(t_j)-\lam_i(t_{j-1}) \right|^2
\eeqN
which implies
\beq\label{d2}
d^2(a,b) \ge \SUM_{j=2}^N d^2(t_j,t_{j-1}).
\eeq
If $d(a,b)\le r$ and $d(t_j,t_{j-1})>\del$ for each~$j$ then~$(N-1)\del^2\le r^2$. Thus
\beqN
\pack(\del,[a,b],d)\le 1 +(r/\del)^2\qt{for }0<\del\le r.
\eeqN
To avoid mess, we simplify the bound to $2(r/\del)^2$.
	\subsection{Chaining bounds}\label{chaining}
In this section we consider a generic
stochastic process $\{X(\theta):\theta\in\Thetabar\}$ whose increments are controlled by the metric~$d$ in the sense that
\beq\label{X.incr}
\PP\{|X(\theta_1)-X(\theta_2)| >  d(\theta_1,\theta_2)x\}
\le C_1e^{-C_2 x}
\qt{for }x\ge 0,
\eeq
for  constants  $C_1$, and~$C_2$.
We establish a one-sided analog of~\cref{Dhat.max} and~\cref{Dmu.max}, 
\bAlignL
\PP&\{\exists \theta\ge \thmu:
 |X(\theta)-X(\theta_\mu)|> 2L(\theta,x,r)\} \le C_1e^{-C_2 x}
\qt{for  $x,r>0$ }
\label{X.weighted}
\\
&\qt{where }L(\theta,x,r) = \left[d^2(\theta,\thmu)+r^2\right]x/r 
\notag
\eAlignL
 We omit the  argument for~$\theta<\thmu$, which is similar.

As explained in Section~\ref{top.level}, we actually only need the inequality for~$r$ equal to~$\max(\Sqrt{3m^*},7x)$, but that
choice plays no role in the derivation of~\cref{X.weighted}.

The method works by cutting the index set into regions where $L(\theta,x,r)$ is approximately constant.
For a given~$r>0$ cover $[\thmu,n]$ by $\cup_{k=1}^m I_k$ where
$I_k=[a_{k-1},a_k]$ and $d^2(a_k,\thmu)=kr^2$ 
for $k=1,\dots,m-1$ and~$d^2(a_m,\thmu)\le kr^2$.
By~\cref{d2}, each~$I_k$ is of $d$-diameter at most~$r$.
Bound the left-hand side of~\cref{X.weighted}  by
\beqN
\SUM_k \PP\{\exists\theta\in I_k: |X(\theta)-X(\theta_\mu)|> 2krx\}.
\eeqN
Here we have used the fact that $d^2(\theta,\thmu)+r^2\ge kr^2$ for all $\theta$ in~$I_k$, with equality at $\theta=a_{k-1}$.
The $k$th term in the sum is less than
\beqN
\PP\{ |X(a_{k-1})-X(\theta_\mu)|> krx\}
+\PP\{\exists\theta\in I_k: |X(\theta)-X(a_{k-1})|> krx\}
\eeqN
By inequality~\cref{X.incr}, the first term is less than~$C_1e^{-C_2\sqrt{k}w}$. The next lemma handles the other contribution.
Taken together they give a bound of the form
$
\SUM_{k\ge 1}C_3\exp(-C_4kx) 
$
for the left-hand side of~\cref{X.weighted}.  If $C_4x\ge 1$ the sum is bounded by a constant times $\exp(-C_4x)$. An increase  in the constant~$C_1$, if necessary, extends the bound to values of~$x$ for which~$C_4x< 1$.

\begin{lemma}\label{stratum}
Suppose  $\{Z(t):t\in T\}$ is a process with continuous sample paths indexed by a set $T$ equipped with a metric~$d$. Suppose also that
\begin{enumerate}
\item 
The diameter of $T$ is $r$ and the packing numbers satisfy
\beqN
\pack(\delta,T,d)\leq C\left(r/\delta\right)^m
\qt{for $0<\del\le r$},
\eeqN
where~$C$ and~$m$ are constants.

\item  
The increments  of $Z$ are controlled by~$d$, in the sense that
\beqN
\PP\{|Z(t_1)-Z(t_2)|>xd(t_1,t_2)\}\leq C_1\exp(-C_2x)
\qt{for all $x\ge 0$}.
\eeqN
\end{enumerate}
Then 
\beqN
\PP\{\SUP_{t\in T}|Z(t)-Z(t_0)|>c_1 x\} \le c_2 e^{-x}
\qt{for all $x\ge 0$}, 
\eeqN
for constants $c_i$ depending on~$C$ and~$m$.
\end{lemma}

\Proof
Define $T_0=\{t_0\}$ and construct packing sets $T_1,T_2,...$ with~
\beqN
N_i=\#T_i\le \pack(\del_i,T,d) \le C2^{mi}
\qt{where~$\del_i= r/2^i$}. 
\eeqN
By construction,
\beqN
\MIN_{t'\in T_i}d(t,t')\leq\delta_i
\qt{for each $t\in T$}.
\eeqN

Let $\{\gamma_i\}_{i\geq 1}$ be a sequence of positive numbers whose value we will later choose. For simplicity of notation write $R_i=\sum_{j\leq i}\gamma_j$ and $R_\infty$ for~$\SUM_{j=1}^\infty\gam_j$. 
Denote $\Delta_i:=\sup_{t_i\in T_i}|Z(t_i)-Z(t_0)|$. By continuity of sample paths,
\beqN
\Del_i \to \Del:=\SUP_{t\in T}|Z(t)-Z(t_0)|
\qt{as $i\to\infty$.}
\eeqN
so that $M_i\to\PP\{\Del> R_\infty\}$. It suffices to 
bound~$M_i:=\PP\{\Delta_i>R_i\}$.


\startDY
Define $\psi_i:T_i\rightarrow T_{i-1}$ as the function that maps $t_i$ to the element in $T_{i-1}$ that is the closest to $t_i$. 
\startDP
Then $\Del_i \le \Del_{i-1}+S_i$ for each~$i$, where
$S_i = \MAX_{t\in T_i}|Z(t_i)-Z(\psi_i t)|$,
which implies the recursive bound
\endDP
\beqN
\PP\{\Delta_i>R_i\}\leq \PP\{\Delta_{i-1}>R_{i-1}\}+\mathbb{P}\{S_i>\gamma_i\}.
\eeqN
Use a union bound to control the second term.
\bAlign
\PP\{S_i>\gamma_i\} 
& \leq \SUM_{t_i\in T_i}\PP\{|Z(t_i)-Z(\psi_it_i)|>\gamma_i\}
\notag\\
& \leq C_1N_i\exp\left(-C_2\gamma_i/\delta_i\right)
\notag\\
& \leq CC_1\exp(im\log 2-C_2\gamma_i2^i/r)
\eAlign

Since we eventually want $\sum_{i\geq 1}\mathbb{P}\{S_i>\gamma_i\}$ to be exponentially small, we choose $\gamma_i$ 
so that $\exp(im\log 2-C_2\gamma_i2^i/r)=\exp(-x)/2^i$, i.e.,
\beqN
\gamma_i=\frac{r}{C_2}2^{-i}(i(m+1)\log 2+x).
\eeqN
This choice of $\gamma_i$ ensures that the tail probability is small enough, but still we do not want $R_i=\sum_{j\leq i}\gamma_j$ to diverge as $i$ grows. Check
\beqN
R_i=\SUM_{j\leq i}\gamma_j=\frac{r}{C_2}\SUM_{j\leq i}[2^{-j}(j(m+1)\log 2+x)] \leq C_3+C_4x.
\eeqN
Here $C_4$ is a universal constant, and $C_3$ only depends on $m$. When $x\geq 1$, we can absorb $C_3$ into the $C_4x$ term. 
\startDP
In summary,
\beqN
M_i=\PP\{\Delta_i> R_i\}\leq\SUM_{j\geq 1}e^{-x}/2^j=e^{-x}.
\eeqN
If  $c_2=e$ then the upper bound~$c_2e^{-x}$ also covers the $0<x<1$ case. Let $i$ go to infinity to complete the proof.
\endDP
\endDY
\endProof

\bibliographystyle{chicago}%
\bibliography{DBP}%

\end{document}